\documentclass[11pt]{article}
\usepackage{a4wide}
\usepackage[all]{xy}
\usepackage[T1]{fontenc}
\usepackage[latin1]{inputenc}
\usepackage[cm]{aeguill}
\usepackage[francais]{babel}
\usepackage{amsmath}
\usepackage{amssymb}
\usepackage{amsthm}
\usepackage{nicefrac}
\usepackage{mathrsfs}
\newtheorem{theoreme}{Théorème}[section]
\newtheorem{proposition}[theoreme]{Proposition}
\newtheorem{corollaire}[theoreme]{Corollaire}
\newtheorem{lemme}[theoreme]{Lemme}

\theoremstyle{definition}
\newtheorem{paragraphe}[theoreme]{}

\newtheorem{remarque}[theoreme]{Remarque}

\renewcommand{\phi}{\varphi}                         
\newcommand{\perv}{\,{}^\mathfrak p\hspace{-0.06cm}}
\DeclareMathOperator{\GL}{GL}
\DeclareMathOperator{\codim}{codim}
\DeclareMathOperator{\End}{End}
\DeclareMathOperator{\Aut}{Aut}

\DeclareMathOperator{\Ker}{Ker}
\renewcommand{\Im}{\,\mathrm{Im}\,}
\DeclareMathOperator{\Coker}{Coker}

\DeclareMathOperator{\Id}{Id}
\DeclareMathOperator{\Supp}{Supp}

\newcommand{\C}{\mathbb{C}}
\renewcommand{\H}{\mathbb{H}}
\renewcommand{\L}{\mathbb{L}}
\newcommand{\N}{\mathbb{N}}
\renewcommand{\P}{\mathbb{P}}
\newcommand{\Q}{\mathbb{Q}}
\newcommand{\R}{\mathbb{R}}
\newcommand{\V}{\mathbb{V}}
\newcommand{\W}{\mathbb{W}}
\newcommand{\Z}{\mathbb{Z}}
\newcommand{\HR}{\mathscr{H}}
\newcommand{\LR}{\mathscr{L}}
\newcommand{\OR}{\mathscr{O}}
\newenvironment{demo}
{
{{\bf Démonstration.}}} {}

\title{Sections hyperplanes à singularités simples et exemples de variations de structure de Hodge}
\author{Damien Mégy\\ \footnotesize{Institut Fourier, 100 rue des maths BP74, 38402 St Martin d'Hères Cedex, France}\\\footnotesize{damien.megy@ujf-grenoble.fr}}
\date{}

\begin{document}
\maketitle
\abstract{\noindent On construit des variétés projectives lisses complexes de dimension $3$ à $6$ munies de variations de structure de Hodge, en généralisant une construction que J. Carlson et C. Simpson avaient effectuée en dimension $2$. Ensuite, on étudie quelques unes de leurs propriétés, en particulier leurs groupes de cohomologie.\\
\\
\begin{center}{\bf Abstract}\end{center}
\vskip 0.5\baselineskip
We construct smooth complex projective varieties of dimension $3$ to $6$ with variations of Hodge structure, by generalizing an example of J. Carlson and C. Simpson in dimension $2$. Then, we study some of their properties, in particular their cohomology groups.\\
\\
\footnotesize{Mots-clés: variations de structure de Hodge, cohomologie d'intersection, modules de Hodge}\\
\\
\footnotesize{Classification mathématique: 14D05, 14D07, 55N33}

}

\section{Introduction}\label{s1}

Les exemples de variations de structure de Hodge (VSH) sur des bases compactes sont relativement rares, bien que l'on sache que ces objets sont fondamentaux dans l'étude des groupes kählériens, depuis les travaux de C. Simpson (\cite{Simpson92}). L'objet de cet article est la construction puis l'étude de VSH sur certaines variétés projectives lisses.

\begin{paragraphe}

\label{notations}
Soit $X$ une variété projective lisse sur $\C$ de dimension \emph{impaire} $n+1\geq 3$ plongée dans un espace projectif $\P$, linéairement normale, et $\P^\vee$ l'espace dual de $\P$, qui paramètre les sections hyperplanes de $X$. Considérons la famille universelle de sections hyperplanes
\[\mathfrak{X} := \{ (x,H) \in X\times\P^\vee | x\in H \}\]
et notons $\pi : \mathfrak{X} \to \P^\vee$ et $a: \mathfrak{X} \to X$ les projections sur chacun des deux facteurs. Si $s \in \P^\vee$, on note $X_s$ la section hyperplane de $X$ correspondant à $s$. Soit $U$ l'ouvert de $\P^\vee$ qui paramètre les sections hyperplanes lisses, et $X^\vee$ son complémentaire, la variété duale.\\

Si $\W$ est une variation de structure de Hodge rationnelle et polarisable sur $X$ (voir §\ref{s2}), alors on note 
\[ \V := \Coker\left(\underline{H^n(X,\W)}_U \to (R^n\pi_*a^*\W)_{|U}\right)\]
le système local de cohomologie évanescente, qui est une VSH sur $U$. Si $\W$ est le système local trivial $\underline{\C}_{X}$, alors la fibre de $\V$ en un point $u \in U$ est simplement la cohomologie évanescente  $H^n(X_u,\C)_{ev}:= \Coker\left(H^n(X,\C) \to H^n(X_u,\C)\right)$ de la section hyperplane.\\
\end{paragraphe}

Dans \cite{Simpson93}, poursuivant une idée de J. Carlson (publié dans \cite{CT93}), C. Simpson considère un plan projectif générique $\P^2 \subset \P^\vee$. Si le plongement $X\subset \P$ est suffisamment ample ($3$-jet ample, voir §\ref{s2}), l'intersection $X^\vee \cap\P^2$ est une courbe irréductible dont les seules singularités sont des n\oe uds et des cusps. Les groupes de monodromie locaux (voir §\ref{s4}) de $\V$ autour des points de cette courbe sont finis, ce qui permet de construire un revêtement ramifié $Y$ de $\P^2$ auquel il est possible d'étendre $\V$ en une VSH $\V_Y$. Sous certaines conditions sur les données initiales $(X,\OR_X(1),\W)$, la VSH $\V_Y$ a des propriétés intéressantes:  un gros groupe de monodromie et une application des périodes génériquement immersive (voir \cite{CT93}, \cite{CT99}). En application de ses résultats sur les images directes supérieures de fibrés harmoniques, Simpson montre également que la VSH ainsi construite est non rigide si l'on choisit bien $X$ et $\W$ (\cite{Simpson93}, th. 9.2).\\

Le premier résultat de cet article est l'extension de cette construction en dimensions $3$ à $6$.\\

\begin{theoreme} --- \label{théorème_construction_Y} Soient $X$ et $\OR_X(1)$ comme en \ref{notations},  $m$ un entier compris entre $3$ et $6$, $\P^m \subset \P^\vee$ un sous-espace projectif générique de dimension $m$, et $U_m:=U\cap \P^m$ l'ouvert paramétrant des sections hyperplanes lisses. La restriction de $\V$ à $U_m$ est encore notée $\V$.

Si $\OR_X(1)$ est $m$-jet ample, il existe une variété quasiprojective lisse $\widetilde{U_m}$, un revêtement galoisien (fini de groupe $G$) $p : \widetilde{U_m}\to U_m$, une compactification lisse $i:\widetilde{U_m}\hookrightarrow Y$, finie sur $\P^m$, tel que $\V_Y:=i_*p^*\V$ soit une VSH sur $Y$.\\
\end{theoreme}

Les variations de structure de Hodge ainsi construites héritent des propriétés démontrées par Carlson et Simpson en dimension deux, concernant leur rigidité, leur groupe de monodromie, ainsi que l'immersivité générique de leur application des périodes. On en déduit des résultats de non factorisation par des variétés de dimension inférieure similaires à ceux de \cite{Simpson93}.

Remarquons que le fait d'étendre la construction de Carlson en dimension $3$ et $5$ permet d'itérer la construction, et d'obtenir beaucoup d'autres exemples de variétés projectives munies de VSH intéressantes.

La preuve du théorème \ref{théorème_construction_Y} consiste à voir que les conditions sur $\OR_X(1)$ et sur $m$ suffisent pour que les groupes de monodromie locaux de $\V$ au voisinage de $\P^m \setminus U$ soient des groupes finis. En effet, par le théorème \ref{ouvert_V_l}, les sections hyperplanes paramétrées par $s \in\P^m$ n'ont alors que des singularités simples (voir §\ref{s2}), et la monodromie locale de $\V$ se déduit des représentations de monodromie de ces singularités. Or, celles-ci ont une image finie.  Ceci permet de construire le revêtement étale. La construction d'une compactification lisse et finie sur $\P^m$ utilise une description locale de la variété duale donnée dans le paragraphe \ref{s3}.\\

On étudie ensuite le couple $(Y,\V_Y)$ par des méthodes cohomologiques. La conjecture de Carlson-Toledo dit que si le groupe fondamental $\pi_1(Z)$ d'une variété kählérienne compacte $Z$ n'est pas fini, alors son second nombre de Betti est virtuellement non nul, c'est-à-dire qu'il existe un sous-groupe $\Gamma \subset \pi_1(Z)$ d'indice fini tel que $H^2(\Gamma,\R) \neq 0$.

Il est bien connu (Reznikov) que si $\V_Z$ est un système local sur $Z$ tel que $H^1(Z,\V_Z)\neq 0$, alors $H^2(\pi_1(Z),\R)\neq 0$.

Reprenons les notations du théorème \ref{théorème_construction_Y}. Dans la suite de l'article, on montre sous certaines hypothèses (voir le théorème \ref{th_cohomologie_invariante}) que la partie $G$-invariante du groupe $H^k(Y,\V_Y)$ est isomorphe au groupe de cohomologie d'intersection $IH^k(\P^m,\V)$, que l'on sait décrire (théorème \ref{SH_IH_sous_famille}). Il s'avère en particulier que $H^1(Y,\V_Y)^G \simeq IH^1(\P^m,\V) \simeq H^{n+1}(X,\W)_{prim}$, ce qui fournit donc un critère de non annulation de $H^2(\pi_1(Y),\R)$. De cette façon, on obtient beaucoup d'exemples de variétés $Y$ de dimension $2$ à $6$ pour lesquelles la conjecture de Carlson-Toledo est vérifiée. Pour des énoncés précis et un peu plus généraux, voir le théorème \ref{th_cohomologie_invariante}.\\

L'article est organisé comme suit. Quelques  définitions et notations sont rappelées au paragraphe \ref{s2}. Au paragraphe \ref{s3}, on démontre quelques propriétés de la variété duale $X^\vee$. Au paragraphe \ref{s4}, on étudie la monodromie locale de $\V$ au voisinage des sections hyperplanes à singularités simples, puis, en utilisant les résultats du paragraphe \ref{s3}, on démontre le théorème \ref{théorème_construction_Y}. Dans le paragraphe \ref{s5}, on étudie la décomposition de Saito de la famille d'hypersurfaces (théorème \ref{decomposition_E_i_explicite}) et on prouve  le théorème \ref{SH_IH_sous_famille} qui décrit la cohomologie d'intersection de $\V$. Dans le paragraphe \ref{s6}, on donne une condition suffisante pour qu'une extension intermédiaire de système local coïncide avec son extension au sens des faisceaux (proposition \ref{prop_extension_intermédiaire}), puis on prouve un théorème  de décomposition pour les familles génériques de petite dimension (théorème \ref{prop_image_directe_décomposable}). Enfin, au paragraphe \ref{s7}, on applique tous ces résultats à l'étude de la cohomologie invariante $H^*(Y,\V_Y)^G$ et on prouve le théorème \ref{th_cohomologie_invariante}.\\

Remarquons que dans le théorème \ref{théorème_construction_Y}, on pourrait remplacer le sous-espace générique $\P^m$ par n'importe quelle sous-variété projective lisse de $\P^\vee$, générique et de dimension $m$. De plus le théorème \ref{théorème_construction_Y} reste valide même si $\W$ est une variation de structure de Hodge polarisée complexe (non nécessairement rationnelle).\\

Cet article est issu de la thèse \cite{MegyThese}, préparée à l'Institut Fourier sous la direction de Philippe Eyssidieux.  Je suis extrêmement reconnaissant à Philippe Eyssidieux, ainsi qu'à Stéphane Guillermou et Chris Peters, pour de nombreuses conversations mathématiques. Je remercie également le rapporteur pour m'avoir aidé à rendre le texte plus lisible.

\section{Notations et rappels}\label{s2}

\begin{paragraphe}\label{s2.vsh}
Un système local en espaces vectoriels complexes $\V$ sur une variété $X$ est (sous-jacent à) une variation de structure de Hodge complexe polarisée de poids $w$ s'il existe une filtration décroissante $F^\bullet$  de $\V \otimes\OR_X$ par sous-fibrés holomorphes et une forme hermitienne $S$ non dégénérée plate telles que
\begin{itemize}
\item[(i)]{le fibré différentiel $V$ sous-jacent à $\V$ se décompose en une somme directe de sous-fibrés $V=\oplus_{p+q=w} V^{p,q}$ vérifiant $F^p = \oplus_{p'\geq p}V^{p',q}$;}
\item[(ii)]{la décomposition de $V$ en somme directe soit $S$-orthogonale, et $(-1)^pS$ soit définie positive sur $V^{p,q}$;}
\item[(iii)]{la transversalité de Griffiths soit vérifiée: $DF^p \subset F^{p-1} \otimes \Omega_X^1$.\\}
\end{itemize}
Une VSH complexe est dite rationnelle si $\V$  a une structure rationnelle plate $\V_\Q \subset \V$ telle que $V^{p,q} = \overline{V^{q,p}}$. Une polarisation rationnelle est une forme $\V_\Q\times\V_\Q\to\V_\Q$ bilinéaire, $(-1)^w$-symétrique, telle que $S(\cdot,\overline\cdot)$ (ou $i.S(\cdot,\overline\cdot)$) soit une polarisation complexe. Dans toute la suite, les VSH seront toujours polarisées.\\
\end{paragraphe}

\begin{paragraphe}\label{compl_loc}
Soit $f \in \C\{x_0, ...x_n\}$  un germe de singularité d'hypersurface. On note $Tju(f)$ son idéal de Tjurina, c'est-à-dire l'idéal engendré par $f$ et ses dérivées partielles, et, si la singularité est isolée, on appelle nombre de Tjurina et on note $\tau$ la codimension (finie) de l'idéal de Tjurina. Comme base d'une déformation miniverselle de la singularité $f$, on peut prendre un voisinage de l'origine dans un supplémentaire $\C^\tau$ de l'idéal de Tjurina. On note alors $\Sigma \subset \C^\tau$ le diagramme de bifurcation de $f$, et $\rho_f : \pi_1(\C^\tau\setminus \Sigma)\to \GL(E)$ la représentation de monodromie de $f$ (voir \cite{AVG}, tome II, §3).

Si de plus $f$ est  une singularité  simple (c'est-à-dire de type ADE par \cite{Arnold}, th. 2.10) et de dimension paire, la représentation de monodromie de $f$ se factorise par le groupe de Coxeter  $W$ de  même type (\cite{LIE}) et par sa représentation comme groupe fini de réflexions (voir \cite{AVG}, tome II). Le sous-groupe distingué $\Ker \rho_f \;\triangleleft\;\pi_1(\C^\tau\setminus \Sigma)$ définit un revêtement étale galoisien $p$ de $\C^\tau\setminus \Sigma$  qui peut alors être complété en un revêtement galoisien ramifié $\widetilde p : \C^\tau \to \C^\tau$. Essentiellement, le $\pi_1$ coïncide avec le groupe de tresses $B(W)$ associé à $W$, le noyau de $\rho_f$ est exactement le groupe des tresses pures, et $\widetilde p$ est le quotient sous l'action de $W$ comme groupe de réflexions, voir par exemple \cite{AVG} ou \cite{MegyThese} p. 35-36.

\end{paragraphe}

\begin{paragraphe}\label{s2.jet_ample}
Enfin, si $\LR$ est un faisceau inversible sur $X$, et $m$ est un entier, on dit que $\LR$ est $m$-jet ample (\cite{BS}) si pour $l$-uplet d'entiers naturels $(m_i)_{1\leq i\leq l}$ tels que $\sum m_i=m+1$ et tout $l$-uplet $(x_i)$ de points distincts de $X$, l'évaluation des multijets
\[ H^0(X,\LR) \to \bigoplus_{i=1}^l \frac{\OR_{X,x_i}}{\mathfrak{m}_{x_i}^{m_i}}\]
est surjective.
\end{paragraphe}

\section{Propriétés de la variété duale}\label{s3}

\begin{theoreme} \label{ouvert_V_l}--- Soit $(X,\OR_X(1))$ comme en \ref{notations}, et $m$ un entier compris entre $2$  et $6$. Supposons que $\OR_X(1)$ soit $\max(3,m)$-jet-ample. Alors, il existe un sous-ensemble fermé $Z_m$ de $X^\vee$ de codimension au moins $m+1$ dans $\P^\vee$ tel que pour tout $v \in X^\vee\setminus Z_m$, les singularités de la section hyperplane $X_v$ soient isolées, simples, et la somme de leurs nombres de Tjurina soit inférieure ou égale à $m$.\end{theoreme}

\begin{demo} --- \'Ecrivons $\P^\vee$ comme l'union disjointe:
\[\P^\vee = \underbrace{X^\vee_{ni} \sqcup X^\vee_{ns} \sqcup X^\vee_{s, \sum \tau_i > m}}_{Z_m} \sqcup V_m, \]
où l'on a noté $X^\vee_{ni}$ l'ensemble des sections hyperplanes avec singularités non isolées, $X^\vee_{ns}$ l'ensemble des sections hyperplanes à singularités isolées mais dont certaines ne sont pas simples, $X^\vee_{s, \sum \tau_i > m}$ l'ensemble des sections hyperplanes n'ayant que des singularités simples et dont la somme des nombres de Tjurina est supérieure à $m+1$, et enfin $V_m$  l'ouvert complémentaire des trois premiers sous-ensembles. Montrons que $Z_m$ est de codimension au moins $m+1$.\\
\begin{itemize}
\item[\emph{\'Etape 1: $\codim X^\vee_{ni} \geq m+1.\label{codimension_non_isolées}$}]{$ $\\
Ce résultat est connu, voir par exemple \cite{Schnell08}, lemme 6.5.2.\\
}
\item[\emph{\'Etape 2: $\codim(X^\vee_{ns}) = 7$.}]{$ $\\ \label{codimension_non_simples}
Il suffit de démontrer que le cône $\C^*.X^\vee_{ns}$ est de codimension $7$ dans $H^0(X,\OR_X(1))$. Soit $\sigma$ une section globale de $\OR_X(1)$ dont les points singuliers $(x_i)_{i\in I}$ sont isolés, et qui a $l$ singularités non simples, aux points $x_1, ..., x_l$. Soient $U_i$ des petites boules de $X$ centrées aux points $x_i$ et disjointes. Soit $B$ une boule de $H^0(X,\OR_X(1))$ centrée en $\sigma$ telle que si $\nu \in B$, alors les singularités de $\nu$ sont toutes isolées et dans les ouverts $U_i$.

\'Ecrivons $B \cap \C^*.X^\vee_{ns} = \bigcup_{i=1}^l NS_i$, où $NS_i$ est l'ensemble des éléments de $B$ ayant une singularité non simple dans $U_i$, et montrons que $NS_i$ est un sous-ensemble de codimension $7$ de $B$.

Soit $j_{x_i} : H^0(X,\OR_X(1)) \to \OR_{X,x_i}$ l'évaluation des jets en $x_i$, et $H$ son image. On note $B'$ l'image de $B$, et $C$ l'image de $NS_i$. Rappelons que $\OR_X(1)$ engendre ses $3$-jets. Par \cite{Arnold}, th. 2.11, $C$ est de codimension $7$ dans $B'$, et comme $j_{x_i}$ est un isomorphisme sur son image, $NS_i$ est de codimension  $7$ dans $B$. Ainsi, $B \cap \C^*.X^\vee_{ns}$ est une union finie de sous-ensembles de codimension $7$ de $B$, donc il est de codimension $7$.\\
}

\item[\emph{\'Etape 3: codimension des strates restantes.}]{$ $\\
Soient $f_1$, ...$f_k$ des germes de singularités isolées (éventuellement les mêmes) et $\tau_1, ..., \tau_k$ leurs nombres de Tjurina, avec $\sum_{i=1}^k \tau_i > m$. On note $S(f_1, .., f_k)$ ou simplement $S$ l'ensemble des sections globales de $\OR_X(1)$ dont les singularités sont exactement les $(f_i)$. Montrons que  
\begin{align}\codim S \geq m+1.\label{controle_strate_avec_somme_tjurina}\end{align}

Soient $a_1, ..., a_k$ des entiers naturels tels que $a_i\leq \tau_i$ et $\sum_{i=1}^k a_i = m+1$. Alors, comme $\OR_X(1)$ est $m$-jet ample, le morphisme d'évaluation des multijets ci-dessous est surjectif
\[H^0(X,\OR_X(1)) \xrightarrow{\oplus j_{x_i}} \bigoplus_{1\leq i\leq k} \OR_{X,x_i} \to \bigoplus_{1\leq i\leq k} \frac{\OR_{X,x_i}}{\mathfrak{m}_{x_i}^{a_i}}.\]

On en déduit que
\[ \codim\left( \bigcap_{i=1}^k j_{x_i}^{-1}(Tju(f_i)+\mathfrak{m}_{x_i}^{a_i})\right) = \sum_{i=1}^k \codim\left(  j_{x_i}^{-1}(Tju(f_i)+\mathfrak{m}_{x_i}^{a_i})\right).\]
Or, puisque $a_i \leq \tau_i$, on a la minoration
\[ \codim (Tju(f_i)+\mathfrak{m}_{x_i}^{a_i}) \geq a_i.\]
En effet, l'idéal $Tju(f_i)$ est de codimension finie $\tau_i$, donc la suite d'entiers 
\[\left(\codim (Tju(f_i)+\mathfrak{m}_{x_i}^{l})\right)_{l\in\N}\]
est strictement croissante jusqu'à atteindre $\tau_i$, puis elle est constante. Finalement,
\[ \codim\left( \bigcap_{i=1}^k j_{x_i}^{-1}(Tju(f_i)+\mathfrak{m}_{x_i}^{a_i})\right) \geq \sum_{i=1}^k a_i=m+1\]

Revenons maintenant à l'ensemble $S$, que l'on suppose non vide. Soit $f \in S$. Au voisinage de $f$, l'ensemble $S$ est l'intersection d'ensembles analytiques $S_i$ formés des sections globales de $\OR_X(1)$ proches de $f$ qui ont exactement une singularité $f_i$ au voisinage de $x_i$. L' espace tangent de Zariski $T_{S_i,f}$ de $S_i$ en $f$ s'envoie dans $Tju(f_i)$ par l'application $j_{x_i}$ de jet en $x_i$ (voir \cite{AVG} ou \cite{Dimca87}, 6.52 p. 92). On en déduit 
\[ \bigcap_{i=1}^k T_{S_i,f} \subset \bigcap_{i=1}^k j_{x_i}^{-1}(Tju(f_i)+\mathfrak{m}_{x_i}^{a_i}).\]
En prenant la codimension des deux termes, on a bien
\[ \codim S \geq m+1. \]
}
\item[\emph{Dernière étape.}]{$ $\\

Il existe une collection finie de singularités simples à laquelle appartiennent toutes les singularités des sections hyperplanes $X_p$ pour tout $p \in X^\vee_{s, \sum \tau_i > m}$. En effet, par \cite{Dimca86}, la somme des nombres de Milnor des singularités de $X_p$ est égale à la multiplicité de $X^\vee$ en $p$. On en déduit qu'il existe $\mu_{max}$ tel que pour tout $p \in X^\vee_s$, les singularités aient toutes un nombre de Milnor inférieur à $\mu_{max}$. Or, il n'y a qu'un nombre fini de telles singularités simples. Ceci montre aussi que le nombre de singularités d'une section hyperplane $X_p$ avec $p \in X^\vee_s$ est borné indépendamment de $p$.

On en déduit que $X^\vee_{s, \sum \tau_i > m}$ est une union finie d'ensembles de la forme de la forme $S(f_1, ... f_k)$, où les germes de singularités simples $f_i$ ont des nombres de Tjurina vérifiant l'inégalité $\sum \tau_i \geq m+1$. Or on a vu que les $S(f_1, ... f_k)$ sont de codimension au moins $m+1$. Finalement, $Z_m$ est bien de codimension au moins $m+1$.\qed\\
}
\end{itemize}
\end{demo}

\begin{paragraphe}\label{structure_locale}
Soit $s \in X^\vee$, tel que les points singuliers de la section hyperplane $X_s$ soient isolés. On les note $x_1$, ..., $x_k$, et on note $f_1, ...,f_k$ les germes de singularités correspondants. On note enfin $\tau_1, ...,\tau_k$ les nombres de Tjurina de ces singularités.

Le germe de $\P^\vee$ au point $s$ est la base d'une déformation de la variété singulière $X_s$. Cette déformation induit une déformation de chaque singularité de $X_s$. Par la propriété universelle des déploiements miniversels de singularités (\cite{Looijenga84}, p. 101), il existe des morphismes de germes d'espaces analytiques de $(\P^\vee,s)$ vers les bases des déploiements miniversels des singularités $f_i$, c'est-à-dire, quitte à choisir une base de $\OR_{X,x_i} / Tju(f_i)$:

\[ (\P^\vee,s) \xrightarrow{g_i} (\C^{\tau_i},0).\]

Les différentielles en $s$ de ces morphismes de germes sont déterminées de manière unique. Ainsi, il existe un morphisme
\[ g : (\P^\vee,s) \to \prod_{i=1}^k (\C^{\tau_i},0)\]
dont la différentielle en $s$ est détérminée de manière unique (et non nulle car $\OR_X(1)$ est très ample).\\

Notons $\Sigma_i \subset \C^{\tau_i}$ les diagrammes de bifurcation des singularités, $U_i$ les ouverts complémentaires, $U=\prod U_i$ et $\Sigma$ le complémentaire de $U$ dans $\prod_{i=1}^k \C^{\tau_i}$, autrement dit
\[\Sigma:= \left\{ (t_i) \in \prod_{i=1}^k \C^{\tau_i}, \exists i, t_i\in\Sigma_i\right\}.\]
\end{paragraphe}

Le résultat qui suit décrit la variété duale au point $s$.\\

\begin{proposition} \label{prop_structure_variete_duale}--- Le germe de $X^\vee$ en $s$ est exactement le germe de $g^{-1}(\Sigma)$.\end{proposition}
\begin{demo} ---
Si $s'$ est assez proche de $s$, les singularités de $X_{s'}$ proviennent toutes par déformation des singularités de $X_s$. Ainsi, les $s'$ proches de $s$ qui appartiennent à la variété duale sont dans l'image réciproque par $g$ d'au moins un diagramme de bifurcation $\Sigma_i$.\qed
\end{demo}

\begin{proposition} \label{prop_lisse}--- Si le faisceau $\OR_X(1)$ est $((\sum\tau_i)-1)$-jet ample, alors le morphisme $g$ est lisse.\end{proposition}
\begin{demo} --- La condition de positivité sur $\OR_X(1)$ implique que  $H^0(X,\OR_X(1)) \to \bigoplus \frac{\OR_{X,x_i}}{\mathfrak{m}_{x_i}^{\tau_i}}$ est surjective, et d'autre part, $\bigoplus \frac{\OR_{X,x_i}}{\mathfrak{m}_{x_i}^{\tau_i}} \to \bigoplus \frac{\OR_{X,x_i}}{Tju(f_i)} $ est également surjective car  $\mathfrak{m}_{x_i}^{\tau_i} \subset Tju(f_i)$. La surjectivité de $dg$ en $s$ en résulte.\qed\\\end{demo}

Des deux propositions précédentes on tire le

\begin{corollaire} \label{cor_pi_1_local} ---
Soit $B_s$  une petite boule centrée en $s$ et $U_s := U \cap B_s$. Sous les hypothèses de \ref{prop_lisse}, le morphisme
\[ g_* : \pi_1(U_s) \to \prod_{i=1}^k \pi_1(\C^{\tau_i} \setminus \Sigma_i)\]
est un isomorphisme.
\end{corollaire}

\section{Généralisation d'une construction de Carlson et Simpson}\label{s4}

Soit $(X,\OR_X(1),\W)$ comme en \ref{notations}, et $u\in U$. La représentation de monodromie de $\V$ est notée $\rho: \pi_1(U,u) \to \GL(\V_u)$. Pour $s\in X^\vee$, on note $B_s$ une petite boule de $\P^\vee$ centrée en $s$ et $U_s := U \cap B_s$. Si $u \in U_s$, on note $\rho_s : \pi_1(U_s,u) \to \GL(\V_u)$ la représentation de monodromie locale de $\V$ au voisinage de $s$, c'est-à-dire la représentations de monodromie de $\V$ restreint à $U_s$. Son image est le groupe de monodromie local.

\begin{proposition} \label{prop_monodromie_locale}--- Soit $s \in X^\vee$ tel que la section hyperplane correspondante $X_s \subset X$  n'ait que des singularités simples. Soit $m$ un entier supérieur à la somme des nombres de Tjurina des singularités de $X_s$. Si $\OR_X(1)$ est $m$-jet ample,  la représentation de monodromie locale de $\V$ au voisinage de $s$ est isomorphe au produit des représentations de monodromie des singularités de $X_s$.
\end{proposition}
 
\begin{demo} --- Donnons la preuve dans le cas de la représentation de monodromie locale sur le groupe d'homologie \[ H_n(X_u,\W)_{ev}:=\Ker(H_n(X_u,\W)\to H_n(X,\W)).\]

Soient $x_1, ...x_l$ les points singuliers de $X_s$, $f_i$ les germes de singularité correspondants et $(B_i)$ une famille de petites boules disjointes centrées en $x_i$. Soit $i$ un entier compris entre $1$ et $l$. Alors $F_i:=X_u \cap B_i$ s'identifie à \og la\fg$\:$ fibre de Milnor de la singularité de $X_s$ en $x_i$. L' inclusion $j_i : F_i \hookrightarrow X_u$ induit une flèche en homologie:
\[{j_i}_*: H_n(F_i,\W) \to H_n(X_u,\W).\]
La flèche composée 
\[ H_n(F_i,\W)\xrightarrow{{j_i}_*} H_n(X_u,\W) \to H_n(X,\W)\]
est nulle; en effet, elle se factorise par $H_n(B_i,\W)$, qui est trivial car la boule $B_i$ est contractile. On en déduit que l'image de ${j_i}_*$ est incluse dans l'homologie évanescente $H_n(X_u,\W)_{ev}$.\\

Les espaces vectoriels $H_n(F_i,\W)$ et $H_n(X_u,\W)_{ev}$ sont munis de formes hermitiennes obtenues en combinant la forme d'intersection sur les cycles topologiques et la polarisation de $\W$ sur les sections de $\W$. Ces formes hermitiennes seront notées $\langle\:,\:\rangle$. La forme hermitienne correspondante sur $H_n(F_i,\W) \simeq H_n(F_i,\C)\otimes \W_{x_i}$ est le produit de la forme d'intersection et de la forme hermitienne $S$ sur la fibre $\W_{x_i}$.

La forme $S$ est par hypothèse non dégénérée sur $\W_{x_i}$. De plus, la forme d'intersection sur l'homologie de la fibre de Milnor $F_i$ est non dégénérée. Alors, la forme hermitienne $\langle\:,\:\rangle$ sur $H_n(F_i,\W)$ est non dégénérée. D'autre part, le morphisme ${j_i}_*$ est compatible aux deux formes hermitiennes ci-dessus; il est donc injectif. 

Soient $i$ et $k$ deux entiers différents compris entre $1$ et $l$, $\alpha \in \Im({j_i}_*)$ et $\beta \in \Im({j_k}_*)$; ce sont des cycles de $X_u$ à valeurs dans $\W$ sont le support est disjoint, puisque $F_i$ et $F_k$ sont disjoints. Donc $\langle\alpha,\beta\rangle=0$. On en déduit que les images  des morphismes ${j_i}_*$ et ${j_k}_*$ sont orthogonales pour la forme $\langle\:,\:\rangle$ sur $H_n(X_u,\W)$. D'autre part ces images sont non dégénérées. Elles sont donc d'intersection triviale. Autrement dit, on a une injection compatible aux structures hermitiennes
\[ \bigoplus_i^{\bot} H_n(F_i,\W) \hookrightarrow H_n(X_u,\W).\]
Dans la suite on écrira par abus de langage les espaces $H_n(F_i,\W)$ comme des sous-espaces de $H_n(X_u,\W)_{ev}$.\\

Revenons à la représentation de monodromie locale $\rho_s : \pi_1(U_s,u) \to \GL(H_n(X_u,\W)_{ev})$. Le groupe $\pi_1(U,u)$ est engendré par des lacets tournant autour des points du lieu lisse de $X^\vee$. Ces lacets agissent sur $H_n(X_u,\W)_{ev}$ par des réflexions de Picard-Lefschetz. Parmi ces lacets, ceux qui sont homotopes à des lacets inclus dans $U_s$ engendrent $\pi_1(U_s,u)$ et agissent par des transformations de Picard-Lefschetz de la forme 
\[ \gamma \cdot \alpha = \alpha \pm \langle \alpha,\delta\rangle\delta,\]
où  $\delta$ est la sphère évanescente correspondant au lacet $\gamma$ (située dans une des fibres de Milnor $F_i$), et l'écriture $\langle\alpha,\delta\rangle$ désigne l'intersection du $\W$-cycle $\alpha$ avec la sphère $\delta$, c'est-à-dire la section de $\W$ sur $\Supp\alpha \cap \delta$ suivante: $\langle \Supp\alpha,\delta\rangle.\alpha_{|\delta}$. Avec ces notations, $\langle \alpha,\delta\rangle\delta$ est bien un $n$-cycle à valeurs dans $\W$.\\

La formule de Picard-Lefschetz montre que $\pi_1(U_s,u)$ stabilise $H_n(F_i,\W)$ dans $H_n(X_u,\W)_{ev}$, ce qui donne des sous-représentations que l'on note 
\[ \rho_{s,i}: \pi_1(U_s,u) \to H_n(F_i,\W)).\]
Par construction, la représentation $\rho_{s,i}$ provient de la représentation de monodromie $\rho_{f_i}$ de la singularité $f_i$, au sens où le diagramme ci-dessous est commutatif (avec les notations de \ref{structure_locale} - \ref{prop_lisse}):
\[\xymatrix{
\pi_1(U_s,u) \ar[r]^{{g_i}_*}\ar[d]^{\rho_{s,i}} &\pi_1(\C^{\tau_i}\setminus \Sigma_i,g_i(u))\ar[d]^{\rho_{f_i}}\\
\GL(H_n(F_i,\C)\otimes \W_{x_i}) & \GL(H_n(F_i,\C))\ar[l]_{-\otimes \Id}
}\]

Par ailleurs la représentation quotient
\[ \rho'_s : \pi_1(U_s,u) \to \GL(H_n(X_u,\W)_{ev} / (\oplus_i H_n(F_i,\W)))\]
est triviale, également par la formule de Picard-Lefschetz. La représentation $\rho_s$ est donc isomorphe au produit des  $\rho_{s,i}$.\\

Enfin, comme $\OR_X(1)$ est $m$-jet ample, le morphisme 
\[g_*=\prod {g_i}_* : \pi_1(U_s,u) \to \prod \pi_1(\C^{\tau_i}\setminus \Sigma_i,g_i(u))\]
est un isomorphisme (corollaire \ref{cor_pi_1_local}). Les représentations $\rho_{s,i}$ sont donc isomorphes aux représentations de monodromie des singularités de $X_s$ en $x_i$.\qed\\
\end{demo}

\begin{paragraphe}\label{revetement}
Soit $\widetilde{U_s}\to U_s$ le revêtement galoisien étale de $U_s$ associé au sous-groupe distingué $\Ker \rho_s$ de $\pi_1(U_s)$, et $p_i: \widetilde{\C^{\tau_i}\setminus \Sigma_i} \to \C^{\tau_i}\setminus \Sigma_i$ celui associé au sous-groupe $\Ker \rho_{f_i}$.

Sous les hypothèses de la proposition \ref{prop_monodromie_locale}, l'isomorphisme $g_*$ identifie donc $\Ker \rho_s$ à $\prod\left(\Ker\rho_{f_i}\right)$. Le revêtement $\widetilde{U_s}\to U_s$ est exactement donné par le produit fibré
\[
\xymatrix{
\widetilde{U_s} \ar[d]\ar[r]\ar@{}[dr]|{\square} & \prod (\widetilde{\C^{\tau_i}\setminus \Sigma_i}) \ar[d]^{\prod_i p_i}\\
U_s \ar[r]^{g} & \prod_i \left(\C^{\tau_i}\setminus \Sigma_i\right)
}\]
\end{paragraphe}

\begin{theoreme} (théorème \ref{théorème_construction_Y} de l'introduction) --- Soit $(X,\OR_X(1),\W)$ comme en \ref{notations}, et $m$ un entier compris entre $3$ et $6$. Soit $\P^m \subset \P^\vee$ un sous-espace projectif générique de dimension $m$ et $U_m:=U\cap \P^m$. La restriction de $\V$ à $U_m$ est encore notée $\V$.\ Si $\OR_X(1)$ est $m$-jet ample, il existe une variété quasiprojective lisse $\widetilde{U_m}$, un revêtement étale galoisien $p : \widetilde{U_m}\to U_m$, une complétion projective lisse $i:\widetilde{U_m}\hookrightarrow Y$, et un revêtement (ramifié) $\hat p: Y\to \P^2$ qui complète le précédent, de telle sorte que $\V_Y:=i_*p^*\V$ soit un système local sur $Y$.\end{theoreme}

\begin{demo} --- On procède en deux étapes: dans un premier temps, on construit le revêtement étale $p : \widetilde{U_m}\to U_m$, puis on montre qu'il existe une variété projective lisse $Y$, finie sur $\P^m$, qui complète ce revêtement. Le principe de la construction du revêtement étale $p : \widetilde{U_m} \to U_m$ est déjà présent dans \cite{CT93}, dans le cas où $m=2$, et repose sur le fait que les monodromies locales de $\V$ sont finies à l'infini (théorème \ref{prop_monodromie_locale}). La démonstration du deuxième point est différente, et utilise la comparaison entre singularités simples et groupes de Coxeter de type ADE (§\ref{s2}).\\

\begin{itemize}
\item[\emph{\'Etape 1: construction du revêtement étale.}]{$ $\\
L'image $\Gamma$ de la représentation de monodromie $\rho: \pi_1(U_m,u) \to \GL(V)$ contient par le théorème de Selberg un sous-groupe distingué $G'$ d'indice fini sans torsion. Soit $G = \rho^{-1}(G')$. Soit $p : \widetilde{U_m} \to U_m$ le revêtement galoisien étale associé sous-groupe distingué $G$ de $\pi_1(U_m,u)$. 

Soit $s \in X^\vee \cap \P^m$, $B_s$ une petite boule de $\P^\vee$ centrée en $s$, $B_{m,s}$ son intersection avec $\P^m$, $U_{m,s} = B_s \cap U_m$,  $\widetilde{U_{m,s}}$ une composante connexe de $p^{-1}(U_{m,s})$ et $\widetilde u$ un point de $\widetilde{U_{m,s}}$ se projetant sur $u$. Par la proposition \ref{prop_monodromie_locale}, la restriction de $\V$ à $U_{m,s}$ a une monodromie finie. Ceci, joint au fait que $G'$ est sans torsion, montre que le revêtement $\widetilde{U_{m,s}} \to U_{m,s}$ est donné par le sous-groupe distingué
\[ \Ker(\rho_s: \pi_1(U_{m,s},u) \to \GL(\V_u)).\]
Le système local $p^*\V$ restreint à l'ouvert $\widetilde{U_{m,s}}$ est donc trivial. On en déduit que pour toute complétion projective $Y$ de $\widetilde{U_m}$, le faisceau $\V_Y:=i_*p^*\V$ est un système local sur $Y$.\\}

\item[\emph{\'Etape 2: construction de la complétion lisse.}]{$ $\\
Soit $Y$ la complétion de Grauert-Remmert du revêtement $\widetilde{U_m}$ (voir \cite{SGA1}, théorème 5.4 p. 254). Alors, $Y$ est une variété normale, finie sur $\P^m$, qui complète le revêtement $\widetilde{U_m} \to U_m$, et c'est la seule à isomorphisme près.

Pour montrer que $Y$ est lisse, il suffit de construire, pour tout $s$ comme ci-dessus, une variété lisse $\widetilde{B_{m,s}}$ et un revêtement ramifié fini $\widetilde{B_{m,s}} \to B_{m,s}$ qui complète le revêtement étale $\widetilde{U_{m,s}} \to U_{m,s}$. En effet, toujours par le théorème de Grauert-Remmert, une telle complétion locale $\widetilde{B_{m,s}}$ est unique à isomorphisme près, et provient par restriction de $Y$.\\

La section hyperplane $X_s$ a $l$ singularités, de type ADE. On note encore $\tau_i$ leurs nombres de Tjurina, $\C^{\tau_i}$ des bases de déploiements miniversels, $\Sigma_i \subset \C^{\tau_i}$ les diagrammes de bifurcation.

Par la propriété universelle des déformations miniverselles, on a un morphisme
\[ g' : (\P^m,s) \to \prod_{i=1}^l(\C^{\tau_i},0)\]
qui se factorise en
\[ (\P^m,s) \hookrightarrow (\P^\vee,s) \xrightarrow{g} \prod_{i=1}^l(\C^{\tau_i},0),\]
où $g$ a été étudié en \ref{structure_locale}. Quitte à rétrécir la boule $B_s$, on peut supposer que le germe $g$ est défini sur $B_{s}$.

Comme $\OR_X(1)$ est supposé $m$-jet ample et $\P^m\subset \P^\vee$ est générique, on a des analogues de \ref{prop_lisse}, \ref{cor_pi_1_local}, et de \ref{revetement}:  $g'$ est submersif, le morphisme induit ${g'}_* : \pi_1(U_{m,s},u) \to \prod\pi_1(\C^{\tau_i}\setminus \Sigma_i)$ est un isomorphisme, et le revêtement $\widetilde{U_{m,s}} \to U_{m,s}$ est obtenu  comme produit fibré:
\[
\xymatrix{
\widetilde{U_{m,s}} \ar[d]\ar[r]\ar@{}[dr]|{\square} & \prod (\widetilde{\C^{\tau_i}\setminus \Sigma_i}) \ar[d]^{p}\\
U_{m,s} \ar[r]^{g'} & \prod \left(\C^{\tau_i}\setminus \Sigma_i\right)
}\]

Par \ref{compl_loc}, chaque revêtement étale $p_i$ de $\C^{\tau_i}\setminus \Sigma_i$ déterminé par le noyau de la représentation de monodromie de $f_i$ peut être complété en un revêtement ramifié $\widetilde{p_i}:\C^{\tau_i} \to \C^{\tau_i}$. On note $\widetilde{p} := \prod_i \widetilde{p_i} : \prod_i \C^{\tau_i}\to \prod_i \C^{\tau_i}$. On définit alors la complétion $\widetilde{B_{m,s}}$ par le produit fibré:
\[
\xymatrix{
 \widetilde{B_{m,s}} \ar[d]\ar[r]\ar@{}[dr]|{\square}& \prod \C^{\tau_i}\ar[d]^{\widetilde p} \\
B_{m,s} \ar[r]^{g'}& \prod \C^{\tau_i}
}\]
Alors, $\widetilde {B_{m,s}}$ est lisse et $\widetilde {B_{m,s}} \to B_{m,s}$ est un revêtement fini qui complète $\widetilde{U_{m,s}} \to U_{m,s}$.\qed\\}
\end{itemize}
\end{demo}

\begin{paragraphe}\label{premieres_proprietes}
Les propriétés suivantes de $\V_Y$ sont connues pour $m=2$ et leur démonstration en dimension $\leq 6$ est la même, mutatis mutandis.
\end{paragraphe}

\begin{itemize}
\item[(i)]{Si $\W$ est sous-jacent à une VSH sur $X$, alors le système local $\V_Y$ est sous-jacent à une variation de structure de Hodge sur $Y$ (\cite{Griffiths70}).}
\item[(ii)]{Si $\W=\C$ et $\OR_X(1)$ est assez ample, alors, l'application des périodes de $\V_Y$ est génériquement de rang maximal. En particulier la représentation de monodromie ne se factorise pas par le groupe fondamental d'une variété de dimension inférieure (\cite{Simpson93}, prop. 8.2, voir aussi \cite{CT93}).}
\item[(iii)]{Si $\OR_X(1)$ est assez ample et $\W=\underline{\C}_X$, le groupe de monodromie de $\V_Y$ est d'indice fini dans le groupe orthogonal $\Aut(H^n(X_u,\Z)_{ev},\langle\:,\:\rangle)$ (une section hyperplane a une singularité $U_{12}$, on applique ensuite \cite{Ebeling84} puis \cite{Beauville86}, voir aussi \cite{CT99}, théorème 9.1).}
\item[(iv)]{Si $\W$  fait partie d'une famille non stationnaire de représentations de $\pi_1(X)$, paramétrée par une variété irréductible $T$, alors $\V_Y$ fait également  partie d'une famille non stationnaire de représentations de $\pi_1(Y)$ (\cite{Simpson93}, lemme 1.5, et théorème 5.1 p. 386). Ceci entraîne que la conjecture de Carlson-Toledo est vraie pour $\pi_1(Y)$ (voir par exemple \cite{Reznikov98}, proposition 8.1).\\}
\end{itemize}

L'extension du deuxième point au cas $\W\neq \underline{\C}_X$ fait l'objet de la note \cite{Megy10}. Le résultat obtenu est le suivant:

\begin{theoreme} --- Si $\W$ est une VSH de poids $w$ avec $W^{w,0}\neq 0$ et $\OR_X(1)$ est assez ample, alors, l'application des périodes de $\V_Y$ est génériquement de rang maximal.\qed\\
\end{theoreme}

On peut interpréter le couple $(Y,\V_Y)$ en termes de systèmes locaux ou de VSH sur des champs de Deligne-Mumford: la VSH $\V$ se prolonge en une VSH $\widetilde \V$ sur le champ de Deligne-Mumford $[Y/G]$, dont l'espace grossier est $\P^m$. En un sens, l'objet principal dans cette construction est $([Y/G],\widetilde \V)$ et non $(Y,\V_Y)$. Nous reviendrons brièvement sur cette interprétation dans les sections suivantes, où apparaissent les groupes de cohomologie $H^l([Y/G],\widetilde \V)$. Cependant nous n'étudieront pas directement $[Y/G]$.

\section{Décomposition de Saito sur des sous-familles génériques}\label{s5}

\subsection{Décomposition de Saito sur la famille complète d'hypersurfaces}

Soit $(X,\OR_X(1),\W)$ comme en \ref{notations}. On note $N$ la dimension de $\P^\vee$ et $d_{\mathfrak{X}}$ celle de $\mathfrak X$. La variété d'incidence ${\mathfrak{X}}$ est lisse et le morphisme $a : \mathfrak{X} \to X$ également, donc le tiré en arrière $\W_{\mathfrak{X}} := a^*\W$ est un système local sur $\mathfrak{X}$ et $\W_{\mathfrak{X}}[d_{\mathfrak{X}}]$ est un faisceau pervers sur $\mathfrak{X}$, sous-jacent à un module de Hodge (polarisable) sur $\mathfrak{X}$, de poids $\omega:=w+n+N$. On note $j : U \hookrightarrow \P^\vee$ l'inclusion dans $\P^\vee$ de l'ouvert $U$, complémentaire dans $\P^\vee$ de la variété duale.\\

Dans ce paragraphe, on applique le théorème de décomposition de Saito (\cite{Saito88}, corollaire 5.4.8 p. 992) à l'image directe dérivée $R\pi_*\W_{\mathfrak X}[d_{\mathfrak{X}}]$. Contrairement à l'article \cite{BFNP}, l'objectif n'est pas d'étudier la famille universelle toute entière, mais surtout certaines sous-familles. Il n'est donc pas nécessaire d'expliciter tous les modules de Hodge irréductibles dans la décomposition de Saito, seulement ceux dont le support est de petite codimension. Ceci permet d'affaiblir les hypothèses sur $\OR_X(1)$.

\begin{theoreme} \label{decomposition_E_i_explicite}--- Soit $(X,\OR_X(1),\W)$ comme en \ref{notations}, et  $m > 1$ un entier. Si $\OR_X(1)$ est $m$-jet ample, alors on a la décomposition de Saito suivante:
\begin{align*}
R\pi_*{\W_{\mathfrak{X}}[d_\mathfrak{X}]} =& \left(\bigoplus_{i=-n}^{-1}\underline{H^{n+i}(X,\W)}[N-i]\right)\\
&\oplus j_{!*}\V[N] \oplus \underline{H^n(X,\W)}[N]\oplus R\\
&\oplus \left(\bigoplus_{i=1}^{n} \underline{H^{n-i}(X,\W)}[N-i]\right),
\end{align*}
où $R$ est un faisceau pervers dont le support est de codimension au moins $m+1$.
\end{theoreme}

\begin{demo} --- Le morphisme $\pi$ est propre de dimension relative $n$. Le théorème de décomposition de Saito  pour l'image directe dérivée de $\W_{\mathfrak X}$ s'écrit :

\begin{align}
\label{decomposition_E_i}R\pi_*{\W_{\mathfrak X}[d_{\mathfrak{X}}]} &\simeq \bigoplus_{i=-n}^{n} E_i[-i],
\end{align}
où les $E_i$ sont des modules de Hodge de poids $\omega+i$, dont les faisceaux pervers sous-jacents sont $\perv\HR^i(R\pi_*{\W_{\mathfrak X}[d_{\mathfrak{X}}]})$.\\
 
\begin{itemize}

\item[\emph{\'Etape 1: calcul des $E_i$ pour $i\neq 0$}]{$ $\\
Si $i < 0$, alors $E_i=\underline{H^{n+i}(X,\W)}[N]$ et $E_{-i}=\underline{H^{n+i}(X,\W)}[N](i)$. Le premier point  découle du théorème de Lefschetz faible relatif, qui résulte de \cite{BBD}, th. 4.1.1 (l'image directe dérivée par un morphisme affine est $t$-exacte à droite). Il est démontré explicitement pour les modules de Hodge dans \cite{BFNP} pour $\W=\underline{\C}_X$. Le second point est une conséquence du premier par l'intermédiaire du théorème de Lefschetz difficile de Saito (\cite{Saito88}, th. 5.3.1 p. 977). Le seul module de Hodge qu'il reste donc à calculer est $E_0$.\\}

\item[\emph{\'Etape 2: les facteurs directs de $E_0$ de codimension $0$ sont $j_{!*}\V[N]$ et $\underline{H^n(X,\W)}[N]$.}]{$ $\\
Restreignons $E_0$ à l'ouvert $U$, complémentaire de la variété duale. Comme $\pi$ est lisse au-dessus de $U$, $E_0|_{U}$ est à décalage près un système local. Soit $u$ un point de $U$, $X_u$ la section hyperplane associée et notons toujours $\W$ la restriction de $\W$ à $X_u$. Alors, on a  $H^n(X_u,\W)\simeq H^n(X_u,\W)_{ev}\oplus H^n(X,\W)$. De là, on tire l'isomorphisme de faisceaux pervers
\[ E_0|_{U} \simeq \V[N]\oplus \underline{H^n(X,\W)}_U[N].\]
Donc, $E_0$ contient les extensions intermédiaires $j_{!*}\left(\V[N]\right)$ et $j_{!*}\left(\underline{H^n(X,\W)}_U[N]\right) = \underline{H^n(X,\W)}[N]$ en facteurs directs, et les autres facteurs directs sont supportés hors de l'ouvert $U$, c'est-à-dire dans la variété duale.\\}

\item[\emph{\'Etape 3: facteurs directs de $E_0$ de codimension $1$.}]{$ $\\
L'entier $n$ étant pair, il n'y a pas de facteur irréductible non trivial de $E_0$ supporté en codimension $1$.

Ce résultat est connu, voir par exemple \cite{BFNP} lorsque $\W=\underline{\C}_X$. La preuve, qui ne nécessite pas de considérer les structures de Hodge portées par les différents objets,  est la suivante. Soit $p$ un point général (lisse) de la variété duale $X^\vee$. Considérons un pinceau de Lefschetz $\P^1 \subset \P^\vee$ passant par $p$. C'est une droite projective qui intersecte transversalement $X^\vee$ en son lieu lisse. Notons $j_1$ l'inclusion $\P^1 \setminus X^\vee \subset \P^1$, et $\pi_1:\mathfrak{X}_{\P^1} \to \P^1$ la restriction de $\pi : \mathfrak X \to \P^\vee$ au-dessus de $\P^1$. Par \cite{SGA7}, XVIII, 6.3, si $n$ est pair, on a la propriété suivante (notée \og (A)\fg$\:$ dans \cite{SGA7} XVIII, 5.3.2): le morphisme d'adjonction
\[R^l{\pi_1}_*\W_{\mathfrak{X}} \longrightarrow {j_1}_*{j_1}^*R^l{\pi_1}_*\W_{\mathfrak{X}}\]
est un isomorphisme de faisceaux, pour tout entier $l$. En particulier, si $l=n+1$, on obtient, en prenant la tige au point $p$,  l'isomorphisme 
\[ H^{n+1}(X_p,\W) \simeq H^{n+1}(X_u,\W)^\gamma,\]
où $u \in \P^1 \setminus X^\vee$ est une valeur régulière de $\pi$ proche de $p$, et $\gamma$ est la monodromie autour de $p$, agissant sur $H^{n+1}(X_u,\W)$. En appliquant le théorème de Lefschetz difficile à coefficients semisimples, on obtient finalement l'isomorphisme d'espaces vectoriels
\begin{align}\label{isom_SGA7} H^{n+1}(X_p,\W) \simeq H^{n-1}(X,\W).\end{align}

D'autre part, par le théorème de changement de base propre, $H^{n+1}(X_p,\W)$ est isomorphe à la tige en $p$ du faisceau $\HR^{-N+1}\left(R\pi_*\W_{\mathfrak X}[d_{\mathfrak{X}}]\right)$. La décomposition \ref{decomposition_E_i} donne alors
\[\HR^{-N+1}\left(R\pi_*\W_{\mathfrak X}[d_{\mathfrak{X}}]\right) \simeq \HR^{-N+1}\left(\bigoplus_{i=-n}^n E_i[-i]\right)=\bigoplus_{i=-n}^n\HR^{-N+1-i}\left( E_i\right).\]
Pour $i$ différent de $0$, on a vu que $E_i$, en tant que complexe de faisceaux, est  un système local constant placé en degré $-N$. Il n'a donc de faisceau de cohomologie qu'en degré $-N$. Il y a par conséquent au plus deux termes non nuls dans la somme directe précédente:
\[\HR^{-N+1}\left(R\pi_*\W_{\mathfrak X}[d_{\mathfrak{X}}]\right) \simeq \HR^{-N}(E_1) \oplus \HR^{-N+1}(E_0).\]
Le premier terme est isomorphe au faisceau constant de tige $H^{n-1}(X,\W)$ en vertu de la première étape. En prenant la tige au point $p$, l'isomorphisme \ref{isom_SGA7} donne donc
\begin{align}\label{nul} \left(\HR^{-N+1}(E_0)\right)_p = 0.\end{align}

Supposons maintenant qu'il existe  un facteur irréductible $E$ non trivial de $E_0$ supporté en codimension un. Son support est donc dense dans la variété duale $X^\vee$. Soit $p$ un point général de son support. Alors, $\left(\HR^{-N+1}(E)\right)_p$ est non trivial, donc $\left(\HR^{-N+1}(E_0)\right)_p$ non plus, ce qui contredit \ref{nul}.\\}

\item[\emph{\'Etape 4: facteurs directs de $E_0$ de codimension $>1$}]{$ $\\
Si $E$ est un module de Hodge irréductible facteur direct de $E_0$ tel que $\codim \Supp E >1$, alors son support est contenu dans le lieu des sections hyperplanes à singularités non isolées.

Ce résultat est dû à Fakhruddin si $\W=\underline{\C}_X$ (voir \cite{BFNP}, remarque 5.12 et théorème A.1). La preuve dans le cas d'une VSH $\W$ non triviale est la même (et c'est essentiellement le même raisonnement qu'à l'étape précédente), à un détail près.

Soit $S$ le support de $E$, $c>1$ la codimension de $S$, et $p$ un point général de $S$. Alors la tige de $E$ en $p$ est un espace vectoriel non nul placé en degré $-N+c$. Calculons $\HR^{-N+c}(R\pi_*\W_{\mathfrak{X}}[d_{\mathfrak{X}}])$ de la même façon qu'à l'étape précédente. Le théorème de décomposition donne, après oubli des structures de Hodge:
\begin{align*}
\HR^{-N+c}\left(R\pi_*\W[d_{\mathfrak{X}}]\right) &\simeq \HR^{-N+c}\left(\bigoplus E_i[-i]\right)\\
&\simeq \HR^{-N+c}(E_0) \oplus \left(\bigoplus_{i \neq 0}\HR^{-N+c-i}(E_i)\right)
\end{align*}
Pour la même raison qu'à l'étape précédente, la somme $\bigoplus_{i \neq 0}\HR^{-N+c-i}(E_i)$ ne contient qu'un seul terme non nul, lorsque $i=c$. On en déduit
\begin{align*}
\HR^{-N+c}\left(R\pi_*\W[d_{\mathfrak{X}}]\right)&\simeq \HR^{-N+c}(E_0) \oplus \HR^{-N}(E_c)\\
&\simeq\HR^{-N+c}(E_0) \oplus \underline{H^{n-c}(X,\W)}.\\
\end{align*}

Ainsi, la tige en $p$ de $\HR^{-N+c}(R\pi_*\W_{\mathfrak{X}}[d_{\mathfrak{X}}])$ contient en facteur direct l'espace vectoriel $H^{n-c}(X,\W)$, ainsi que la tige de $E$ en $p$, qui est un espace vectoriel non trivial. D'autre part, par le changement de base propre pour l'inclusion $p \hookrightarrow \P^\vee$,  la tige en $p$ de $\HR^{-N+c}(R\pi_*\W_{\mathfrak{X}}[d_{\mathfrak{X}}])$ est isomorphe à $H^{n+c}(X_p,\W)$. On en déduit que $H^{n+c}(X_p,\W)$ n'est pas isomorphe à $H^{n-c}(X,\W)$.\\

Ceci n'est possible que si $X_p$ est à singularités non isolées. En effet, on a le lemme suivant:
\begin{lemme} --- Si $X_p$ est une section hyperplane irréductible à singularités isolées de $X$ et $c$ un entier strictement supérieur à $1$, alors $H^{n+c}(X_p,\W) \simeq  H^{n-c}(X_p,\W)$.\qed\\
\end{lemme}

Ce résultat vient remplacer le théorème de dualité de Kaup invoqué par Fakhruddin. Il est certainement bien connu, au moins dans le cas $\W = \underline\C_X$. La preuve dans le cas général en est une simple adaptation; on a notamment les isomorphismes suivants:
\begin{itemize}
\item{si $k<n$, alors $IH^k(X_p,\W)\simeq H^k(X_p^{lisse},\W)$;}
\item{si $k>n$, alors $IH^k(X_p,\W)\simeq H_c^k(X_p^{lisse},\W)\simeq H^k(X_p,\W)$;}
\item{si $k<n-1$, alors $H^k(X_p^{lisse},\W)\simeq H^k(X_p,\W)$.}
\end{itemize}
Les deux premiers points sont dus à Goresky et MacPherson \cite{GM80} dans le cas de l'homologie d'intersection à coefficients constants (voir aussi l'introduction de \cite{Durfee93}, et \cite{Dimca04} pour une démonstration utilisant le langage des faisceaux pervers). Le troisième point remonte également aux travaux de Goresky et MacPherson, et résulte, après excision, de la $(n-2)$ connexité du \emph{link} d'une singularité isolée d'hypersurface (Milnor). Finalement, pour $1<c<n$, on a la chaîne d'isomorphismes
\[H^{n-c}(X_p,\W) \simeq IH^{n-c}(X_p,\W) \simeq IH^{n+c}(X_p,\W) \simeq H^{n+c}(X_p,\W),\]
où le second isomorphisme est obtenu par le théorème de Lefschetz difficile de Saito sur la cohomologie d'intersection à coefficients dans une VSH.\\}

\item[\emph{Fin de la démonstration}]{$ $\\
Supposons qu'il y ait un facteur direct $E$ non trivial dans $E_0$, supporté en codimension strictement positive. Le support d'un tel facteur direct est inclus dans le lieu des sections hyperplanes à singularités non isolées. Comme $\OR_X(1)$ est $m$-jet ample, ce lieu est de codimension au moins $m+1$ (voir l'étape 1 de la démonstration du théorème \ref{ouvert_V_l}).\qed}
\end{itemize}
\end{demo}

\begin{remarque} \label{rq_R=0} --- Si de plus $\OR_X(1)$ est $n$-jet ample, alors le théorème \ref{decomposition_E_i_explicite} implique $R=0$. La raison en est la suivante. Soit 
\[
R\pi_*{\W_{\mathfrak X}[d_{\mathfrak{X}}]} \simeq \bigoplus_{i=-n}^{n} E_i[-i]
\]
la décomposition de Saito. Si un module de Hodge irréductible non trivial $E$ est facteur direct dans un des $E_i$, alors par \cite{Ngo08} théorème 2 p. 188, on a $\codim\Supp E \leq n-|i|$ (Ngo  attribue l'argument à Goresky et MacPherson). En particulier, si $R$ n'est pas trivial, son support doit être de codimension inférieure à $n$. Mais ceci est impossible d'après le théorème. Cette remarque, dans le cas $\W=\underline{\C}_X$ et avec une condition plus forte sur $\OR_X(1)$, est exactement la remarque 5.12 de \cite{BFNP}.\\\end{remarque}

\subsection{Calculs de cohomologie d'intersection}

Dans ce paragraphe, on munit la cohomologie d'intersection $IH^k(\P^m,\V)$ d'une filtration et on calcule ses gradués en fonction de différentes structures de Hodge sur la cohomologie de $X$. Le résultat principal est la preuve du théorème \ref{SH_IH_sous_famille}. On commence par démontrer plusieurs résultats intermédiaires.\\

Le lemme élémentaire suivant permettra de munir la cohomologie d'intersection de $\V$ de filtrations et de calculer les gradués.

\begin{lemme} \label{lemme_filtrations}--- Soit $E$ et $H$ deux espaces vectoriels de dimension finie, et $F^\bullet E$, $F^\bullet H$ des filtrations décroissantes indéxées par $\N$, birégulières \footnote{c'est-à-dire vérifiant $\bigcup_i F^iE = E$ et $\bigcap_i F^i E = \{0\}$, $\bigcup_i F^iH = H$ et $\bigcap_i F^i H = \{0\}$}. Soit $r : E\to H$ un morphisme compatible aux filtrations. On note $r_i : F^iE \to F^i H$ et $q_i : Gr^iE\to Gr^iH$ les morphismes induits.

Supposons que $r$ soit injectif, et qu'il existe $i_0$ tel que $q_i$ soit bijectif si $i>i_0$, injectif si $i=i_0$, et surjectif si $i<i_0$. Alors $\Coker q_{i_0}$ porte une filtration dont les gradués sont $\Coker r$ et $\Ker q_i$, pour $0\leq j\leq i_0-1$.\end{lemme}

\begin{remarque} --- Le même énoncé est vrai dans la catégorie des structures de Hodge.\end{remarque}

\begin{demo} du lemme \ref{lemme_filtrations} --- On considère le morphisme de suites exactes courtes:
\[\xymatrix{
0 \ar[r] & F^{i+1} H \ar[r] & F^{i} H \ar[r] & Gr^{i} H \ar[r] & 0\\
0 \ar[r] & F^{i+1} E \ar[r]\ar[u]^{r_{i+1}} & F^{i} E \ar[r]\ar[u]^{r_i} & Gr^{i} E \ar[r]\ar[u]^{q_i} & 0\\
}\]
Les flèches $r_i$ sont toutes injectives, donc le lemme du serpent donne une suite exacte:
\[ 0\to \Ker q_i \to \Coker r_{i+1} \to \Coker r_i \to \Coker q_i \to 0.\]
Si $i>i_0$, $q_i$ est par hypothèse un isomorphisme et donc $\Coker r_{i+1} \simeq \Coker r_i$. Pour $i$ assez grand, on a $\Coker r_i=0$ donc par récurrence descendante, on obtient:
\[ \forall i>i_0, \: \Coker r_i=0.\]
Si $i=i_0$, on obtient 
\[\Coker r_{i_0} \simeq \Coker q_{i_0}.\]
Si $i<i_0$, $q_i$ est surjectif et on a la suite exacte courte:
\[ 0 \to \Ker q_i \to \Coker r_{i+1} \to \Coker r_i \to 0.\]
On définit alors une filtration de $K:=\Coker r_{i_0}$ par $F^{-1}K:=K$ et:
\[ F^j K := \Ker(\Coker r_{i_0}\xrightarrow{\Id}\Coker r_{i_0} \to \Coker r_{i_0-1}  \to ... \to \Coker r_j), \: 0\leq j\leq i_0,\]
(en particulier, $F^{i_0}K=\{0\}$). Alors, on a $Gr^{-1}K \simeq \Coker r_0$ et, pour $0 \leq j\leq i_0-1$, 
\begin{align*}
Gr^jK &= \frac{\Ker(\Coker r_{i_0} \to ... \to \Coker r_j)}{\Ker(\Coker r_{i_0} \to ... \to \Coker r_{j+1})}\\
&\simeq \Ker(\Coker r_{j+1} \to \Coker r_j)\\
&\simeq \Ker q_j.\qed\\
\end{align*}
\end{demo}

\begin{lemme} \label{lemme_decomposition_1}--- L'image directe dérivée $R{a_m}_*{a_m}^*\W$ est décomposable, et on a, dans la catégorie dérivée $D_c^b(X)$, la décomposition
\[ R{a_m}_*{a_m}^*\W \simeq \left(\bigoplus_{0\leq i\leq 2m-2} \underline{H^i(\P^m,\C)}\otimes\W[-i] \right) \oplus \underline{H^{2m}(\P^m,\C)}\otimes\W_{|B}[-2m].\]
\end{lemme}

\begin{demo} --- La fibre du morphisme ${a_m}$ au-dessus de $x\in X$ est :
\begin{itemize}
\item{l'espace projectif ${\P^m}$ si $x\in B$;}
\item{un hyperplan projectif de ${\P^m}$ sinon.}
\end{itemize}

Considérons la restriction de ${a_m} : \mathfrak{X}_{\P^m} \to X$ au-dessus de $X\setminus B$. C'est un sous-fibré projectif du fibré trivial $X\times \P^m$. En particulier, ${a_m} : {a_m}^{-1}(X\setminus B) \to X\setminus B$ est projectif lisse. D'après le critère de dégénérescence de Deligne (\cite{Deligne68}), on a une décomposition

\[ R{a_m}_*{a_m}^*\W \simeq \bigoplus_{i=0}^{2m-2} (R^i{a_m}_*{a_m}^*\W)[-i]\]
Dans $D^b(X\setminus B)$. Calculons ces images directes supérieures. Pour $0\leq i\leq 2m-2$ et $x \in X\setminus B$, l'inclusion ${a_m}^{-1}(x) \hookrightarrow \P^m$ induit des isomorphismes en cohomologie:
\[ H^i({a_m}^{-1}(x),{a_m}^*\W) \xleftarrow{\sim} H^i(\P^m,\C)\otimes \W_x.\]
Finalement, on a, sur $X\setminus B$, l'isomorphisme $R^i{a_m}_*{a_m}^*\W\simeq \underline{H^i(\P^m,\C)}\otimes\W$ pour $0\leq i\leq 2m-2$. Si $B$ est vide (c'est-à-dire si $m>n$), ceci termine la démonstration.

Sinon, le théorème de décomposition pour ${a_m} : \mathfrak{X}_{\P^m} \to X$ dit que l'image directe dérivée se décompose en somme directe de faisceaux pervers sur $X$ décalés convenablement. La restriction à $X\setminus B$ de cette décomposition est exactement celle obtenue par le critère de dégénérescence Deligne, écrite ci-dessus. On en déduit que les complexes $\underline{H^i(\P^m,\C)}\otimes\W$ pour $0\leq i\leq 2m-2$, sont facteurs directs de $R{a_m}_*{a_m}^*\W$, et que les autres termes apparaissant dans la décomposition de Beilinson-Bernstein-Deligne pour ${a_m}_*{a_m}^*\W$ sont supportés dans $B$.\\

Considérons donc maintenant la restriction de ${a_m} : \mathfrak{X}_{\P^m} \to X$ au-dessus de $B$. C'est le fibré trivial $B\times\P^m$. On en déduit la décomposition
\[ R{a_m}_*{a_m}^*\W = \bigoplus_{i=0}^{2m} (R^i{a_m}_*{a_m}^*\W)[-i]=\bigoplus_{0\leq i\leq 2m} \underline{H^i(\P^m,\C)}\otimes\W_{|B}[-i].\]
Finalement, l'image directe dérivée $R{a_m}_*\W$ sur $X$ admet la décomposition suivante:
\[ R{a_m}_*{a_m}^*\W \simeq \left(\bigoplus_{0\leq i\leq 2m-2, i\text{ pair}} \underline{H^i(\P^m,\C)}\otimes\W[-i] \right) \oplus \underline{H^{2m}(\P^m,\C)}\otimes\W_{|B}[-2m].\qed\]
\end{demo}

\begin{corollaire} \label{corollaire_restriction}--- Soit $l$ un entier positif strictement inférieur à $m$. Alors, le morphisme de restriction
\[H^{n+l}(\P^m\times X,pr_2^*\W) \to H^{n+l}(\mathfrak{X}_{\P^m},{a_m}^*\W)\]
est un isomorphisme, et le morphisme
\[H^{n+m}(\P^m\times X,pr_2^*\W) \to H^{n+m}(\mathfrak{X}_{\P^m},{a_m}^*\W)\]
est injectif de conoyau isomorphe à $H^{n-m}(B,\W)_{ev}:=\Coker (H^{n-m}(X,\W) \to H^{n-m}(B,\W))$ si $m\leq n$, bijectif sinon.
\end{corollaire}

\begin{demo} --- On considère la décomposition du lemme \ref{lemme_decomposition_1}, ainsi que la décomposition pour la projection $pr_2: \P^m\times X \to X$, qui est un fibré projectif:
\[R{pr_2}_*pr_2^*\W  \simeq \left(\bigoplus_{0\leq i\leq 2m} \underline{H^i(\P^m,\C)}\otimes\W[-i] \right).\]
Le morphisme de restriction en cohomologie est compatible avec ces deux décompositions.

Si $m>n$, le lieu de base est vide, et donc pour tout $l\leq m$, le morphisme de restriction
\[H^{n+l}(\P^m\times X,pr_2^*\W) \to H^{n+l}(\mathfrak{X}_{\P^m},{a_m}^*\W)\]
est un isomorphisme.\\

Si $m \leq n$, le lieu de base est non vide et pour tout $l$ on a
\[
\Ker(H^{n+l}(\P^m\times X,pr_2^*\W) \to H^{n+l}(\mathfrak{X}_{\P^m},{a_m}^*\W))\]
\[=\Ker (H^{n+l-2m}(X,\W) \to H^{n+l-2m}(B,\W)),\]
et
\[ \Coker(H^{n+l}(\P^m\times X,pr_2^*\W) \to H^{n+l}(\mathfrak{X}_{\P^m},{a_m}^*\W))\]
\[ =\Coker (H^{n+l-2m}(X,\W) \to H^{n+l-2m}(B,\W)).\]

Par le théorème de Lefschetz faible pour l'inclusion de  $B$ dans $X$, on en déduit que $H^{n+l}(\P^m\times X,pr_2^*\W) \to H^{n+l}(\mathfrak{X}_{\P^m},{a_m}^*\W)$ est un isomorphisme pour $l<m$ et est injectif pour $l=m$.\qed\\
\end{demo}

\begin{theoreme} \label{SH_IH_sous_famille}--- Soit $(X,\OR_X(1),\W)$ comme en \ref{notations}. Soit $m$ un entier compris entre $2$ et $6$. Supposons que $\OR_X(1)$ soit $m$-jet ample. Soit $\P^m$ un sous-espace projectif générique de $\P^\vee$, $B \subset X$ son lieu de base et $X \xleftarrow{a_m} \mathfrak{X}_{\P^m} \xrightarrow{\pi_m} \P^m$ la famille correspondante de sections hyperplanes. La restriction de $\V$ à $U\cap\P^m$ est toujours notée $\V$.

Soit $l$ un entier. Si $l<m$, ou si $l=m$ et $m>n$, la structure de Hodge sur  $IH^l(\P^m,\V)$ porte une filtration dont les gradués successifs sont isomorphes à 
\[ H^i(\P^m)\otimes H^{n+i-l+2}(X,\W)_{prim}(i-l+1),\]
pour $i$ variant entre $0$ et $l-1$. Si $l=m$ et que $m\leq n$, la filtration a un cran de plus et le gradué supplémentaire est isomorphe à $H^{n-m}(B,\W)_{ev}(-m)$.\end{theoreme}

\begin{demo} --- On cherche à appliquer le lemme \ref{lemme_filtrations}. Notons $pr_2 : \P^m\times X \to X$ la deuxième projection et $a_m : \mathfrak{X}_{\P^m} \to X$ sa restriction à $\mathfrak{X}_{\P^m}$.\\

Les morphismes $\pi:\mathfrak{X}_{\P^m} \to \P^m$ et $pr_1 : \P^m\times X \to \P^m$ donnent des suites spectrales de Leray perverses
\[ \perv E_2^{p,q}(\pi) \Rightarrow H^{p+q}(\mathfrak{X}_{\P^m},{a_m}^*\W),\]
\[ \perv E_2^{p,q}(pr_1) \Rightarrow H^{p+q}(\P^m\times X,pr_2^*\W),\]
qui dégénèrent en $E_2$.\\

Le morphisme $pr_1 : \P^m\times X \to \P^m$ est projectif lisse, donc la suite spectrale de Leray perverse coïncide avec la suite spectrale de Leray classique et on a :
\[ \perv E_2^{p,q}(pr_1) \simeq E_2^{p,q}(pr_1) = H^p(\P^m,\C)\otimes H^q(X,\W).\]

D'autre part, la décomposition du théorème \ref{decomposition_E_i_explicite} peut être restreinte de $\P^\vee$ à $\P^m$, et on a:
\[ \perv E_2^{p,q}(\pi) \simeq \H^{p-N}(\P^m,(E_{q-n})_{|\P^m}).\]

Utilisons maintenant la fonctorialité de la suite spectrale de Leray perverse. On note $E=H^{n+l}(\P^m\times X,pr_2^*\W)$, $H=H^{n+l}(\mathfrak{X}_{\P^m},a_m^*\W)$, $r:E \to H$ le morphisme de restriction,  $F^\bullet E$ et $F^\bullet H$ les filtrations de Leray perverses sur  $E$ et $H$. Ce sont des filtrations décroissantes indexées par l'indice $i$, avec $0\leq i\leq n+l$. La suite de Leray perverse est fonctorielle en structures de Hodge, ce qui implique que $r$ est compatible aux filtrations de Leray perverses sur $E$ et sur $H$. Ceci se traduit par des diagrammes commutatifs à lignes exactes:
\[\xymatrix{
0 \ar[r] & F^{i+1} H \ar[r] & F^{i} H \ar[r] & Gr^{i} H \ar[r] & 0\\
0 \ar[r] & F^{i+1} E \ar[r]\ar[u]^{r_{i+1}} & F^{i} E \ar[r]\ar[u]^{r_i} & Gr^{i} E \ar[r]\ar[u]^{q_i} & 0\\
}\]

Pour appliquer le lemme \ref{lemme_filtrations}, il reste à vérifier les hypothèses sur les morphismes $r$ et $q_i$.\\

Par le corollaire \ref{corollaire_restriction}, on sait que $r$ est bijectif si $l<m$ ou $m>n$, et qu'il est injectif, de conoyau isomorphe à $H^{n-m}(B,\W)_{ev}$, si $l=m$ et que $m\leq n$.\\

Vérifions maintenant les hypothèses sur les morphismes 
\[q_i : Gr^iH^{n+l}(\P^m\times X,pr_2^*\W) \to Gr^{i}H^{n+l}(\mathfrak{X}_{\P^m},a_m^*\W).\]
 Les gradués sont calculables explicitement car les suites dégénèrent en $E_2$:
\[ Gr^{i}H \simeq \perv E_2^{i,n+l-i}(\pi)\simeq \H^{i-m}(\P^\vee,(E_{l-i})_{|\P^m}[m-N]),\text{ et}\]
\[ Gr^{i}E \simeq E_2^{i,n+l-i}(pr_1)\simeq H^i(\P^\vee,\C)\otimes H^{n+l-i}(X,\W).\]

Si $l-i<0$, on a $(E_{l-i})_{|\P^m}[m-N]\simeq \underline{H^{n+l-i}(X,\W)}[m]$ et donc:
\[ Gr^{i}H \simeq \perv E_2^{i,n+l-i}(\pi) \simeq H^i(\P^m,\C)\otimes H^{n+l-i}(X,\W).\]

Comme $\P^m$ est générique, il ne rencontre pas le support du module de Hodge $R$ (voir le théorème \ref{decomposition_E_i_explicite}), et donc $(E_0)_{|\P^m}[m-N]\simeq \underline{H^n(X,\W)}[m] \oplus (j_{!*}\V[N])_{|\P^m}[m-N]$. De plus, l'inclusion $i :\P^m\hookrightarrow \P^\vee$ est transverse à $X^\vee$, donc le foncteur $i^*[N-m]$ envoie $j_{!*}\V[N])$ sur un faisceau pervers simple. On en déduit l'isomorphisme
\begin{align}
(j_{!*}\V[N])_{|\P^m}[m-N] \simeq {j_m}_{!*}(\V_{|U_m}[m]).\label{changement_base_transverse}
\end{align}

On en déduit que
 
\[ Gr^{l}H \simeq \perv E_2^{i,n}(\pi)  \simeq H^l(\P^m,\C)\otimes (H^{n}(X,\W)\oplus  IH^l(\P^m,\V).\]
Enfin, si $k-i>0$, on a :
\[ Gr^{i}H \simeq \perv E_2^{i,n+l-i}(\pi)  \simeq H^i(\P^m,\C)\otimes H^{n-l+i}(X,\W)(i-l).\]

Les noyaux et conoyaux de $q_i$ sont alors donnés par les théorèmes de Lefschetz. Pour $l+n\geq i>l$, le morphisme 
$q_i$ est un isomorphisme. Pour $i=l$, $q_i$ est injectif, et on a $\Coker q_l = IH^l(\P^m,\V)$. Enfin, pour $0 \leq i <l$, $q_i$ est surjectif et son noyau est isomorphe à 
\[ H^i(\P^m)\otimes H^{n+i-l+2}(X,\W)_{prim}(i-l+1).\]

Le lemme \ref{lemme_filtrations} s'applique et donne le résultat.\qed\\
\end{demo}

\begin{remarque} --- Pour terminer cette section, énonçons sans démonstration deux résultats similaires au théorème \ref{SH_IH_sous_famille}. Les notations sont inchangées.
\begin{itemize}
\item[a) ]{Supposons que $\OR_X(1)$ soit $n$-jet ample, et soit $l<n$ un entier positif. La structure de Hodge sur $IH^l(\P^\vee,\V)$ porte une filtration par sous-structures de Hodge dont les gradués successifs sont isomorphes à \[H^i(\P^\vee)\otimes H^{n+i-l+2}(X,\W)_{prim}(i-l+1),\]
pour $i$ variant entre $0$ et $l-1$.\qed\\}

\item[b) ]{Soit $l$ un entier strictement inférieur à $n$. Si $\W=\underline{\C}_X$ et $\OR_X(1)$ est assez ample pour que le théorème de connexité de  Nori (\cite{Nori}) s'applique, la structure de Hodge mixte (SHM) sur $H^l(U,\V)$ porte une filtration de SHM dont les gradués sont isomorphes aux structures mixtes 
\[H^i(U)\otimes H^{n+i-l+2}(X,\C)_{prim}(i-l+1),\] 
pour $i$ variant entre $0$ et $l-1$.\qed\\}
\end{itemize}

La démonstration du point a) est semblable à celle du théorème \ref{SH_IH_sous_famille}, et utilise le lemme \ref{lemme_filtrations}, la décomposition de Saito pour la famille complète d'hypersurfaces lorsque $R=0$ (ce qui est le cas si $\OR_X(1)$ est $n$-jet ample, voir remarque \ref{rq_R=0}), et la dégénérescence de la suite spectrale  de Leray pour $a : \mathfrak{X}\to X$ qui est une fibration lisse en espaces projectifs. Ce résultat figure dans \cite{BFNP} pour $l=1$  lorsque $\W=\underline{\C}_X$. La démonstration du point b) utilise le lemme \ref{lemme_filtrations} et le théorème de connexité de Nori. Le résultat est d'ailleurs présent dans \cite{Nori} pour $l=1$. Pour plus de détails sur les preuves de ces deux énoncés, voir \cite{MegyThese}, p. 71 et p. 78.\end{remarque}

\section{Décomposition}\label{s6}

Rappelons pour commencer le résultat fondamental suivant (\cite{SGA7}, XVIII, 5.6.7): si $(X,\OR_X(1))$ est comme en \ref{notations}, $\P^1\subset \P^\vee$ un pinceau de Lefschetz de sections hyperplanes de $X$, et $\pi_1 : \mathfrak{X}_{\P^1} \to \P^1$ le morphisme projectif dont les fibres sont précisément les sections hyperplanes du pinceau, alors $R{\pi_1}_*\C$ est décomposable dans $D^b_c(\P^1)$; en particulier, la suite spectrale de Leray du morphisme $\pi_1$ dégénère en $E_2$.\\

Dans ce paragraphe, nous généralisons ce théorème aux familles génériques de dimension $2$ à $6$ de section hyperplanes de $X$ (théorème \ref{prop_image_directe_décomposable}).\\

La proposition suivante, aux notations indépendantes, donne une condition suffisante pour que l'extension intermédiaire d'un système local sur certaines variétés coïncide avec son extension au sens des faisceaux. On démontre ensuite le théorème \ref{prop_image_directe_décomposable}.\\

\begin{proposition} \label{prop_extension_intermédiaire}--- Soit $X$ une variété algébrique lisse sur $\C$ quelconque de dimension $d$, $j : U \to X$ l'inclusion d'un ouvert de Zariski, $L$ un système local sur $U$. Supposons qu'il existe $\pi : \widetilde X \to X$ un revêtement galoisien fini de groupe $G$, tel que le revêtement restreint $\pi_U : \widetilde U \to U$ soit étale, et que $L_{\widetilde U}:={\pi_U}^*L$ s'étende un système local $L_{\widetilde X}$ sur $\widetilde X$.\\

Alors, l'extension intermédiaire $j_{!*}(L[d])$ est un faisceau (à décalage près):
\[  j_{!*}(L[d])\simeq (j_{*}L)[d].\]
\end{proposition}

\begin{demo} --- Les morphismes $\pi$ et $\pi_U$ sont finis, donc les foncteurs image directe associés sont exacts sur les faisceaux, égaux à leurs foncteurs dérivés. On les note simplement $\pi_*$ et ${\pi_U}_*$. Notons $\widetilde j$ l'inclusion de $\widetilde U$ dans $\widetilde X$.\\

Le faisceau image directe ${\pi_U}_*L_{\widetilde U}$ est muni d'une action naturelle de $G$, et on a 
\[ L=({\pi_U}_*L_{\widetilde U})^G.\]
On en déduit que $j_{!*}(L[d]) = j_{!*}\left(({\pi_U}_*L_{\widetilde U})^G[d]\right)$, c'est-à-dire, puisque $j_{!*}$ est un foncteur entre faisceaux pervers équivariants sur $U$ et sur $X$ (\cite{BernsteinLunts}, 5.2 p. 41):
\[ j_{!*}(L[d])=j_{!*}\left({\pi_U}_*L_{\widetilde U}[d]\right)^G.\]

Montrons maintenant que l'extension intermédiaire $j_{!*}({\pi_U}_*L_{\widetilde U}[d])$ est isomorphe à $(\pi_*L_{\widetilde X})[d]$. Le morphisme $\pi:\widetilde X\to X$ est fini, donc le foncteur (dérivé) $\pi_*$ est exact pour les $t$-structures perverses sur $D^b_c(\widetilde X)$ et $D^b_c(X)$. Ceci signifie exactement que l'image d'un faisceau pervers est un faisceau pervers. Le faisceau $L_{\widetilde X}$ est un système local sur la variété lisse $\widetilde X$, ce qui fait de $L_{\widetilde X}[d]$ un faisceau pervers sur $\widetilde X$. Par ce qui précède, $\pi_*L_{\widetilde X}[d]$ est un faisceau pervers sur $X$. Ce faisceau pervers est semi-simple (par le théorème de décomposition), et il contient forcément comme sous-objet l'extension intermédiaire de sa restriction à $U$. On a donc un isomorphisme de faisceaux pervers:
\[ \pi_*L_{\widetilde X}[d] \simeq j_{!*}({\pi_U}_*L_{\widetilde U}[d]) \oplus R,\]
où le reste $R$ est un faisceau pervers supporté sur le complémentaire de $U$, qui est une sous-variété de dimension strictement inférieure à $d$. Or, on sait que $(\pi_*L_{\widetilde X})[d]$ est un faisceau placé en degré $-d$. On en déduit que $R=0$ et donc 
\[ j_{!*}({\pi_U}_*L_{\widetilde U}[d]) = \pi_*L_{\widetilde X}[d].\]

Mais alors,  on a : 
\begin{align*}
j_{!*}(L[d])[-d]&=j_{!*}({\pi_U}_*L_{\widetilde U}[d])^G[-d]\\
&=(\pi_*L_{\widetilde X})^G\\
&=\left((\pi_*\circ \widetilde j_*)L_{\widetilde U}\right)^G\\
&=\left((j_*\circ {\pi_U}_*)L_{\widetilde U}\right)^G\\
&= j_*({\pi_U}_*L_{\widetilde U})^G\\
&=j_*L.\qed
\end{align*}
\end{demo}

\begin{theoreme} \label{prop_image_directe_décomposable} --- Soit $(X,\OR_X(1),\W)$ comme en \ref{notations} et $m$ un entier positif compris entre $2$ et $6$.\\
Supposons que $\OR_X(1)$ soit $\max(3,m)$-jet ample. Soit $\P^m$ un sous-espace projectif générique de $\P^\vee$, et $X\xleftarrow{a_m} \mathfrak{X}_m \xrightarrow{\pi_m} \P^m$ la famille de sections hyperplanes paramétrée par $\P^m$. Alors, l'image directe dérivée $R{\pi_m}_*(a_m^*\W)$ est décomposable dans $D^b_c(\P^m,\C)$.\end{theoreme}

\begin{demo} --- Fixons les notations suivantes pour le changement de base :
\[\xymatrix{
\mathfrak{X}_{\P^m} \ar[d]^{{\pi_m}} \ar[r]^{i'} & \mathfrak{X}\ar[d]^{\pi}\\
{\P^m} \ar[r]^{i}& \P^\vee\\
U_m\ar[u]_{j_m} \ar[r]^k &U\ar[u]^{j}
}\]

Par changement de base propre, on a l'isomorphisme de foncteurs:
\[ i^*R\pi_* \simeq R{\pi_m}_*i'^*,\]
donc on a :
\[ i^*R\pi_*\W_{\mathfrak{X}}[n+m] \simeq R{\pi_m}_*i'^*\W_{\mathfrak{X}}[n+m],\]
En écrivant la décomposition de Saito de $R\pi_*\W[N+n]$, on obtient :
\[ R{\pi_m}_*i'^*\W_{\mathfrak{X}}[n+m] \simeq i^*\left( \bigoplus_{i\in\Z} E_i[-i]\right)[m-N],\]
c'est-à-dire, en explicitant la décomposition (voir \ref{decomposition_E_i_explicite}) :
\begin{align}\label{restriction_decomposition_E_i_explicite}
R{\pi_m}_*(\W_{\mathfrak{X}}[m+n]) =& \left(\bigoplus_{i=-n}^{-1}\underline{H^{n+i}(X,\W)}[m-i]\right)\\
&\oplus i^*(j_{!*}\V[N])[m-N] \oplus \underline{H^n(X,\W)}[m]\nonumber\\
&\oplus \left(\bigoplus_{i=1}^{n} \underline{H^{n-i}(X,\W)}[m-i]\right).\nonumber
\end{align}
Par \ref{changement_base_transverse}, on a un isomorphisme de faisceaux pervers sur $\P^m$:
\[i^*(j_{!*}\V[N])[m-N] \simeq {j_m}_{!*}(\V_{|U_m}[m]).\]
Enfin, on rappelle qu'il existe une variété projective lisse $Y$ et un morphisme fini $p : Y \to \P^m$ satisfaisant les conditions suivantes (voir le théorème \ref{théorème_construction_Y}) :
\begin{itemize}
\item{le morphisme $p$ au-dessus de $U_m$ est un revêtement galoisien étale $p : \widetilde{U_m}\to U_m$;}
\item{le tiré en arrière $p_{|\widetilde{U_m}}^*\V$ sur $\widetilde{U_m}$ s'étend en un système local $\V_Y$ sur $Y$.}
\end{itemize}
Par la proposition \ref{prop_extension_intermédiaire}, le faisceau pervers $j'_{!*}(\V_{|U_m}[m])$ est isomorphe au complexe composé du faisceau $j'_{*}\V_{|U_m}$ placé en degré $-m$. Mais alors, la décomposition \ref{restriction_decomposition_E_i_explicite} montre que tous les termes apparaissant dans la décomposition sont des faisceaux, à décalage près. On en déduit que l'image directe $R{\pi_m}_*i'^*\W_{\mathfrak{X}}[n+m]$ est décomposable.\qed\\
\end{demo}

\section{Cohomologie invariante sur $Y$ et applications}\label{s7}

Dans cette section, on calcule la partie $G$-invariante de la cohomologie de $Y$ à coefficients dans le système local $\V_Y$. On donne ensuite deux applications de la non-annulation de certains de ces groupes de cohomologie.

\begin{theoreme} \label{th_cohomologie_invariante}--- Soit $(X,\OR_X(1),\W)$ comme en \ref{notations} et $m$ un entier positif compris entre $2$ et $6$. Supposons que $\OR_X(1)$ soit $\max(3,m)$-jet ample. Soit $\P^m$ un sous-espace projectif générique de $\P^\vee$. Si $l<m$, ou si $l=m$ et $m>n$, la structure de Hodge sur  $H^l(Y,\V_Y)^G$ porte une filtration dont les gradués successifs sont isomorphes à 
\[ H^i(\P^\vee)\otimes H^{n+i-l+2}(X,\W)_{prim}(i-l+1),\]
pour $i$ variant entre $0$ et $l-1$. Si $l=m$ et que $m\leq n$, la filtration a un cran de plus et le gradué supplémentaire est isomorphe à $H^{n-m}(B,\W)_{ev}(-m)$.\end{theoreme}

\begin{demo}  --- 
On a la chaîne d'isomorphismes
\[ H^k(Y,\V_Y)^G \simeq  H^k(\P^m,p_*\V_Y)^G \simeq H^k(\P^m,(p_*\V_Y)^G)\simeq H^k(\P^m,j_*\V).\]
D'après la proposition \ref{prop_extension_intermédiaire}, $j_*\V[m]=j_{!*}\V[m]$, et donc
\[ H^k(Y,\V_Y)^G \simeq IH^k(\P^m,\V).\]
Or, par le théorème \ref{SH_IH_sous_famille}, la cohomologie d'intersection porte une telle filtration.\qed\\
\end{demo}

\begin{paragraphe}
La partie $G$-invariante est sans doute la plus intéressante, car elle reflète uniquement des propriétés de la VSH $\V$ sur $\P^m$ et ne dépend pas du choix du groupe $G$. Cette partie invariante est la cohomologie du champ $[Y/G]$ à coefficients dans $\widetilde \V$:
\[ H^k([Y/G],\widetilde \V) \simeq H^k(Y,\V_Y)^G.\]
En effet, d'une part, la cohomologie du champ quotient $H^k([Y/G],\widetilde\V)$ est exactement la cohomologie équivariante $H^k_G(Y,\V_Y)$, et d'autre part, on peut calculer cette cohomologie équivariante à l'aide de la suite spectrale de Leray équivariante associée  à  l'application $G$-équivariante $Y \to \text{\{pt\}}$ (voir \cite{McCleary01}, p. 500):
\[ E_2^{p,q} = H^p(G,H^q(Y,\V_Y)) \Rightarrow H_G^{p+q}(Y,\V_Y).\]
Le groupe $G$ étant fini, ses groupes de cohomologie supérieurs à coefficients dans une représentation complexe sont nuls, et donc $E_2^{p,q} = 0, \:\forall p>0$. Ceci implique que la suite spectrale dégénère en $E_2$, et donc que
\[ H^k_G(Y,\V_Y) = H^0(G,H^k(Y,\V_Y)) = H^k(Y,\V_Y)^G.\]
\end{paragraphe}

Terminons par deux applications. Si $H^{n+1}(X,\W)_{prim}\neq 0$, alors sous les hypothèses des théorèmes \ref{th_cohomologie_invariante} et \ref{SH_IH_sous_famille}, on obtient $H^1(Y,\V_Y)\neq 0$.\\

La première application est la construction d'un système local sur $Y$ ayant des déformations analytiques. Ceci découle de la proposition suivante, aux notations indépendantes de ce qui précéde.

\begin{proposition} \label{prop_deformations_33} --- Soit $Y$ une variété projective lisse et $\V_Y$ un système local semisimple sur $Y$. Supposons que $H^1(Y,\V_Y)$ soit non nul. Alors, le système local $\L_Y:=\C_Y\oplus \V_Y$ admet des déformations infinitésimales non obstruées, c'est-à-dire qui se relèvent en des déformations analytiques.\end{proposition}

\begin{demo} --- Remarquons que $\End(\L_Y) \simeq \L_Y^\vee\otimes \L_Y$. En utilisant $\L_Y=\C_Y\oplus \V_Y$, on voit que $\V_Y$ est un sous-système local de $\End(\L_Y)$, et donc que $H^1(Y,\V_Y)$ est un sous-espace de $H^1(Y,\End(\L_Y))$. Par hypothèse, on a $H^1(Y,\V_Y)\neq \{ 0\}$, donc $H^1(Y,\End(\L_Y))\neq \{0\}$. Rappelons que l'espace tangent de Zariski à la variété des représentations (ou des sytèmes locaux) est l'espace des cocycles $Z^1(Y,\End(\L_Y))$ (voir \cite{LuMa}, chapitre 2). Par ce qui précède, il est non nul, ce qui signifie exactement que $\L_Y$ admet des déformations infinitésimales.\\

\'Etudions les obstructions à ce que ces déformations infinitésimales se relèvent en des déformations analytiques.\\

L'obstruction de Goldman-Millson est l'application
\[\Phi : H^1(Y,\End(\L_Y)) \to H^2(Y,\End(\L_Y)),\]
donnée sur un $\pi_1(Y)$-cocycle tordu par la formule
\[ u \mapsto \Phi(u)(g_1,g_2) = [u(g_1),g_2\cdot u(g_2)].\]
Le système local $\L_Y$ est semisimple, donc par \cite{Simpson92} et par  \cite{GoldmanMillson}, corollaire du théorème 2 p. 44, une déformation infinitésimale $u \in H^1(Y,\L_Y)$ se relève en une déformation analytique de $\L_Y$ si et seulement si $\Phi(u)=0$.\\

L'obstruction de Goldman-Millson est nulle sur $H^1(Y,\V_Y)$ vu comme sous-espace de $H^1(Y,\End(\L_Y))$. En effet, soit $u \in H^1(Y,\V_Y) \subset H^1(Y,\End(\L_Y))$. Alors le $2$-cocyle $\Phi(u)$, évalué sur un couple $(g_1, g_2)$ d'éléments de $\pi_1(Y,y)$, est donné par le crochet de Lie de deux éléments dans $(\V_Y)_y \subset (\End(\L_Y))_y$.

Or la décomposition $\L_Y = \C\oplus \V_Y$ donne une décomposition par blocs des éléments de $(\End \L_Y)_y$, et les éléments dans $(\V_Y)_y \subset (\End\L_Y)_y$ sont de la forme $\left(\begin{array}{cc}
0&0\\
*&0
\end{array}\right)$. Le crochet de Lie de deux tels éléments est donc nul.\qed\end{demo}

\begin{remarque} Malheureusement cette méthode de permet pas d'obtenir un système local non stationnaire au sens de Simpson (déformations dans l'espace des modules, voir \cite{Simpson93}, p. 346).\\\end{remarque}

Enfin, une autre application concerne la conjecture de Carlson-Toledo sur la variété $Y$:

\begin{proposition} --- Si $H^1(Y,\V_Y)$ est non trivial, alors  $H^2(\pi_1(Y),\C)$ non plus.\end{proposition}

Le principe général de la preuve est connu, on pourra se référer à \cite{Reznikov98}.\\

\begin{remarque} Si $\V_Y$ est non-rigide, alors on sait déjà que la conjecture de Carlson-Toledo est vraie pour $Y$. Par les théorèmes de Simpson, c'est le cas si $\W$ est non-stationnaire sur $X$. Le résultat de la proposition ci-dessus n'utilise pas cette hypothèse.\\\end{remarque}

Pour conclure, remarquons qu'il serait très intéressant d'en savoir plus sur l'éventuelle rigidité de $\V_Y$ dans le cas où $\W$ est rigide.

\medskip

\medskip

Damien Mégy

Institut Fourier, UMR 5582

100, rue des maths BP74

38402 Saint Martin d'Hères Cedex

France

\medskip

damien.megy@ujf-grenoble.fr

\end{document}